%&amstex
\input amsppt.sty
\magnification=\magstep1
\hsize=30truecc
\baselineskip=16truept
\vsize=22.2truecm
\nologo
\TagsOnRight
\pageno=1

\def\Z{\Bbb Z}

\def\Q{\Bbb Q}

\def\C{\Bbb C}

\def\l{\left}
\def\r{\right}
\def\bg{\bigg}
\def\({\bg(}
\def\[{\bg[}
\def\){\bg)}
\def\]{\bg]}
\def\t{\text}
\def\f{\frac}
\def\mo{\roman{mod}}
\def\per{\roman{per}}

\def\se {\subseteq}

\def\sm{\setminus}

\def\bi{\binom}
\def\eq{\equiv}

\def\ls{\leqslant}
\def\gs{\geqslant}
\def\al{\alpha}
\def\ve{\varepsilon}
\def\da{\delta}

\def\Proof{\noindent{\it Proof}}

\def\Remark{\medskip\noindent{\it  Remark}}
\def\Ack{\medskip\noindent {\bf Acknowledgments}}
\topmatter
\hbox {Final version for J. Number Theory, {\tt arXiv:math.CO/0410529}.}
\medskip
\title On various restricted sumsets\endtitle
\author Zhi-Wei Sun$^{1,*}$ and Yeong-Nan Yeh$^2$\endauthor
\leftheadtext{Zhi-Wei Sun and Yeong-Nan Yeh}
\affil
$^1$Department of Mathematics and Institute of Mathematical Science
\\Nanjing University, Nanjing 210093, P. R. China
\\zwsun\@nju.edu.cn
\\{\tt http://pweb.nju.edu.cn/zwsun}
\medskip
$^2$Institute of Mathematics, Academia Sinica,
Taipei, Taiwan
\\mayeh\@ccvax.sinica.edu.tw
\endaffil

\abstract For finite subsets $A_1,\ldots,A_n$ of a field, their
sumset is given by
 $\{a_1+\cdots+a_n:\, a_1\in A_1,\ldots,a_n\in A_n\}$.
 In this paper we study various restricted sumsets
 of $A_1,\ldots,A_n$ with  restrictions of the following forms:
 $$ a_i-a_j\not\in S_{ij},\ \t{or}\ \al_ia_i\not=\al_ja_j,\ \t{or}
 \ a_i+b_i\not\eq a_j+b_j\ (\mo\ m_{ij}).$$
 Furthermore, we gain an insight into relations among recent results
 on this area obtained in quite different ways.
\endabstract
\thanks 2000 {\it Mathematics Subject Classification}.
Primary 11B75; Secondary 05A05, 11C08.
\newline\indent *This author is responsible for communications,
and supported by the National Science Fund for
Distinguished Young Scholars (No. 10425103) and the Key Program of
NSF (No. 10331020) in China.
\endthanks
\endtopmatter
\document
\hsize=30truecc
\baselineskip=16truept

\heading{1. Introduction}\endheading

The additive order of the identity of a field $F$ is either infinite
or a prime, we call it the {\it characteristic} of $F$.

Let $F$ be a field of characteristic $p$, and let $A_1,\ldots,A_n$
be finite subsets of $F$ with $0<k_1=|A_1|\ls\cdots\ls k_n=|A_n|$.
Concerning various restricted sumsets of $A_1,\ldots,A_n$,
the following results are known:

(i) (The Cauchy-Davenport theorem (see, e.g. [N]))
 $$|\{a_1+\cdots+a_n:\, a_1\in A_1,\ldots,a_n\in A_n\}|
 \gs\min\{p,k_1+\cdots+k_n-n+1\}.$$

(ii) (Dias da Silva and Hamidoune [DH]) If $A_1=\cdots=A_n=A$, then
$$|\{a_1+\cdots+a_n:\, a_i\in A,\ a_1,\ldots,a_n\ \t{are distinct}\}|
\gs \min\{p, n|A|-n^2+1\}.$$

(iii) (Alon, Nathanson and Ruzsa [ANR2]) If $k_1<\cdots<k_n$,
then
$$|\{a_1+\cdots+a_n:\,  a_i\in A_i,\ a_i\not=a_j\ \t{if}\ i\not=j\}|
\gs\min\bg\{p,\sum_{i=1}^nk_i-\f{n(n+1)}2+1\bg\}.$$

(iv) (Hou and Sun [HS]) Let $S_{ij}\ (1\ls i,j\ls n,\ i\not=j)$
be finite subsets of $F$ with cardinality $m$.
If $k_1=\cdots=k_n=k$ and $p>\max\{ln,mn\}$ where $l=k-1-m(n-1)$, then
$$|\{a_1+\cdots+a_n:\,  a_i\in A_i,
\ a_i-a_j\not\in S_{ij}\ \t{if}\ i\not=j\}|\gs ln+1.$$

(v) (Liu and Sun [LS]) Let $P_1(x),\ldots,P_n(x)\in F[x]$ be monic and of degree $m>0$.
If $k_n>m(n-1)$, $k_{i+1}-k_i\in\{0,1\}$ for all $i=1,\ldots,n-1$,
 and $p>K=(k_n-1)n-(m+1)\bi n2$,
then we have
$$|\{a_1+\cdots+a_n:\,  a_i\in A_i,\ P_i(a_i)\not=
P_j(a_j)\ \t{if}\ i\not=j\}|\gs K+1.$$

(vi) (Sun [Su]) Let $P_1(x),\ldots,P_n(x)\in F[x]$ have degree $m>0$ with
the permanent of the matrix $(b_j^{i-1})_{1\ls i,j\ls n}$ nonzero,
 where $b_j$ is the leading coefficient of $P_j(x)$.
If $k_1=\cdots=k_n=k>m(n-1)$ and $K=(k-1)n-(m+1)\bi n2<p$, then
$$|\{a_1+\cdots+a_n:\,  a_i\in A_i,\ a_i\not=a_j\ \&\
P_i(a_i)\not=P_j(a_j)\ \t{if}\ i\not=j\}|\gs K+1.$$

While result (ii) was deduced by a deep tool from the representation theory of symmetric groups,
results (iii)--(vi) were obtained from the following basic principle
arising from Alon and Tarsi [AT].

\proclaim{Combinatorial Nullstellensatz {\rm ([A1, A3])}} Let $A_1,\ldots,A_n$ be finite subsets
of a field $F$ with $|A_i|>k_i$ for $i=1,\ldots,n$
where $k_1,\ldots,k_n$ are nonnegative integers.
 If the coefficient
of the monomial $x_1^{k_1}\cdots x_n^{k_n}$ in $f(x_1,\ldots,x_n)\in F[x_1,\ldots,x_n]$
is nonzero and $k_1+\cdots+k_n$ is
the total degree of $f$,
then there are $a_1\in A_1,\ldots,a_n\in A_n$ such that
$f(a_1,\ldots,a_n)\not=0$.
\endproclaim

Lower bounds for various restricted sumsets are usually yielded
with help of the following lemma (or Proposition 2.1 of [HS])
implied by the Combinatorial Nullstellensatz.

\proclaim{Lemma 1.1 {\rm (Alon et al. [ANR1, ANR2])}}
Let $A_1,\ldots,A_n$ be finite nonempty subsets of a field $F$ with $k_i=|A_i|$
for $i=1,\ldots,n$. Let $P(x_1,\ldots,x_n)\in F[x_1,\ldots,x_n]\setminus\{0\}$
and $\deg P\ls\sum_{i=1}^n (k_i-1)$.
If the coefficient of the monomial $x_1^{k_1-1} \cdots x_n^{k_n-1}$
in the polynomial
 $$P(x_1,\ldots,x_n)(x_1+\cdots+x_n)^{\sum_{i=1}^n(k_i-1)-\deg P}$$
does not vanish, then we have
$$|\{a_1+\cdots+a_n:\,  a_i\in A_i, \ P(a_1,\ldots,a_n)\not=0\}|
\gs \sum_{i=1}^n(k_i-1)-\deg P+1.$$
\endproclaim

In the next section, we will develop a general technique
to compute certain coefficients of some polynomials.
Using Lemma 1.1 and our work in Section 2, we will prove the following main
theorems in Section 3.

\proclaim{Theorem 1.1}
Let $F$ be a field of characteristic $p$, and
let $A_1,\ldots,A_n$ be finite nonempty subsets of $F$
with $|A_n|=k$ and $|A_{i+1}|=|A_i|+1$ for $i=1,\ldots,n-1$.
Let $m$ be a positive integer, and let $S_{ij}\se F$ and $|S_{ij}|<2m$
for all $1\ls i<j\ls n$.
If $p>\max\{mn, (k-1)n-mn(n-1)\}$, then we have
$$|\{a_1+\cdots+a_n:\,  a_i\in A_i,\ \t{and}\ a_i-a_j\not\in
S_{ij}\ \t{if}\ i<j\}|\gs (k-1-m(n-1))n+1.$$
\endproclaim

\Remark\ 1.1. Theorem 1.1 can be viewed as a partial generalization
of result (iii).

\proclaim{Theorem 1.2}  Let $k$ and $m$ be positive integers.
Let $A_1,\ldots,A_n$ be subsets of the complex field $\C$ with cardinality $k$, and
let $S_{ij}\ (1\ls i<j\ls n)$ be subsets of $\C$ with at most $2m-1$ elements.
If $\zeta_1,\ldots,\zeta_n$ are distinct $q$th
roots of unity where $q$ is a positive odd integer, then
$$\aligned&\bg|\bg\{\sum_{i=1}^na_i:\,  a_i\in A_i,
\ a_i\zeta_i\not=a_j\zeta_j\ \t{and}
\ a_i-a_j\not\in S_{ij}\ \t{if}\ i<j\bg\}\bg|
\\&\qquad\gs(k-1-m(n-1))n+1.\endaligned$$
\endproclaim

\Remark\ 1.2. A conjecture of Snevily [S] states that
for any cyclic group with odd order
if $A$ and $B=\{b_1,\ldots,b_n\}$ are its subsets with
cardinality $n$ then there is a numbering $\{a_i\}_{i=1}^n$
of the elements of $A$
such that $a_1b_1,\ldots,a_nb_n$ are pairwise distinct.
Using the Combinatorial Nullstellensatz
Alon [A1] confirmed this for the cyclic group $\Z/p\Z$
where $p$ is an odd prime.
Since we can identify a cyclic group of order $q$
with the multiplicative group of all the $q$th roots of unity,
Snevily's conjecture follows from
Theorem 1.2 in the case $k=n$, $m=1$, $A_1=\cdots=A_n=A$ and
 $S_{ij}=\{0\}\ (1\ls i<j\ls k)$, which was first obtained by
Dasgupta, K\'arolyi, Serra and Szegedy [DKSS] in 2001.
Another extension of Snevily's conjecture appeared in [Su].

\proclaim{Theorem 1.3} Let $\al_1,\ldots,\al_n$ be positive reals,
and let $b_1,\ldots,b_n$ be integers.
Let $A_1,\ldots,A_n$ be finite subsets of $\Z$ with cardinality $k>0$.
For $1\ls i<j\ls n$ let $m_{ij}$ be an integer greater than
$2\max\{|x_i-x_j|:\,  x_i\in A_i,\ x_j\in A_j\}|$. Then
 the restricted sumset
$$\bg\{\sum_{i=1}^na_i:\,  a_i\in A_i,\ a_i\al_i\not=a_j\al_j\
\t{and}\ a_i+b_i\not\eq a_j+b_j\ (\mo\ m_{ij})\ \t{if}\ i<j\bg\}$$
has more than $(k-n)n$ elements.
\endproclaim

\proclaim{Corollary 1.1 {\rm (K\'ezdy and Snevily [KS])}}
Let $m$ and $n$ be positive integers with
$n\ls(m+1)/2$. Then, for any $b_1,\ldots,b_n\in\Z$, there exists a
permutation $\sigma$ on $\{1,\ldots,n\}$ such that
$1+b_{\sigma(1)},\ldots,n+b_{\sigma(n)}$ are pairwise distinct modulo
$m$.
\endproclaim
\Proof.  Observe that
$m/2>n-1=\max_{1\ls i\ls j\ls n}(j-i)$.
Applying Theorem 1.3 with $\al_1=\cdots=\al_n=1$ and $A_1=\cdots=A_n=\{1,\ldots,n\}$,
we find that there exists a permutation $\tau$ on $\{1,\ldots,n\}$
such that $\tau(1)+b_1,\ldots,\tau(n)+b_n$ are pairwise distinct
modulo $m$. So the desired result follows. \qed

\Remark\ 1.3. In [KS] Corollary 1.1 was applied to tree embeddings.
Snevily [S] even conjectured that the condition
$n\ls(m+1)/2$ in Corollary 1.1 can be weakened by $n<m$.
\medskip

Let $G=\{a_1,\ldots,a_n\}$ be an additive
abelian group of order $n$,
and let $b_1,\ldots,b_n$ be elements of $G$ with $b_1+\cdots+b_n=0$.
In 1952 M. Hall [H] proved that
there exists a permutation $\sigma$ on $\{1,\ldots,n\}$ such that
$a_1+b_{\sigma(1)},\ldots,a_n+b_{\sigma(n)}$ are pairwise distinct.
Let $\sigma$ be a permutation on $\{1,\ldots,n\}$ such that
$b_1-a_{\sigma(1)},\ldots,b_n-a_{\sigma(n)}$ are pairwise distinct.
Assume that $n>1$ and $a_n=b_n=0$.
Then there exists a permutation $\sigma'$ on $\{1,\ldots,n-1\}$
such that
$$a_{\sigma'(i)}=a_{\sigma(i)}-a_{\sigma(n)}\not=0
\quad \ \t{for every}\ i=1,\ldots,n-1.$$
Since $\{b_i-a_{\sigma'(i)}:\,i=1,\ldots,n-1\}=G\sm\{0\}$,
there is a permutation $\tau$ on $\{1,\ldots,n-1\}$
such that for any $i=1,\ldots,n-1$ we have
$b_i-a_{\sigma'(i)}=a_{\tau(i)}$ and hence
$b_i=a_{\sigma'(i)}+a_{\tau(i)}$. In the case $G=\Z/n\Z$, this
provides a positive answer to an open question of Parker
(cf. [G]).

\heading{2. Relations among coefficients of certain polynomials}
\endheading

 As usual we let $(x)_0=1$ and $(x)_n=x(x-1)\cdots(x-n+1)$ for
 $n=1,2,3,\ldots$. For a polynomial
$$P(x_1,\ldots,x_n)=\sum_{j_1,\ldots,j_n}a_{j_1,\ldots,j_n}x_1^{j_1}\cdots
x_n^{j_n}$$
over a commutative ring,
we write $[x_1^{j_1}\cdots x_n^{j_n}]P(x_1,\ldots,x_n)$
to denote the coefficient $a_{j_1,\ldots,j_n}$.

\proclaim{Lemma 2.1}
 Let  $$P(x_1,\ldots,x_n)=\sum\Sb j_1,\ldots,j_n\gs0\\j_1+\cdots+j_n=m\endSb c_{j_1,\ldots,j_n}
 x_1^{j_1}\cdots x_n^{j_n}\in\C[x_1,\ldots,x_n]$$
 and
 $$P^*(x_1,\ldots,x_n)=\sum\Sb j_1,\ldots,j_n\gs0\\j_1+\cdots+j_n=m\endSb c_{j_1,\ldots,j_n}
 (x_1)_{j_1}\cdots (x_n)_{j_n}.$$
 Suppose that $0\ls\deg P\ls k_1+\cdots+k_n$
 where $k_1,\ldots,k_n$ are nonnegative integers. Then
 $$[x_1^{k_1}\cdots x_n^{k_n}]P(x_1,\ldots,x_n)(x_1+\cdots+x_n)^{k_1+\cdots+k_n-\deg P}$$
coincides with
$$\f{(\sum_{i=1}^n k_i-\deg P)!}{k_1!\cdots k_n!}P^*(k_1,\ldots,k_n).$$
 \endproclaim

 \Proof. Let $K=k_1+\cdots+k_n-\deg P$. Then
  $$\align &[x_1^{k_1}\cdots x_n^{k_n}]P(x_1,\ldots,x_n)(x_1+\cdots+x_n)^K
 \\=&[x_1^{k_1}\cdots x_n^{k_n}]\sum\Sb j_1,\ldots,j_n\gs0\\j_1+\cdots+j_n=m\endSb
 c_{j_1,\ldots,j_n}x_1^{j_1}\cdots x_n^{j_n}
 \sum\Sb i_1,\ldots,i_n\gs0\\ i_1+\cdots+i_n=K\endSb
 K!\f{x_1^{i_1}\cdots x_n^{i_n}}{i_1!\cdots i_n!}
 \\=&K!\sum\Sb j_1,\ldots,j_n\gs0\\j_1+\cdots+j_n=m\endSb
 c_{j_1,\ldots,j_n}\f{(k_1)_{j_1}\cdots(k_n)_{j_n}}{k_1!\cdots k_n!}
 \\=&\f{K!}{k_1!\cdots k_n!}P^*(k_1,\ldots,k_n).
 \endalign$$
 This concludes the proof.
 \qed

Let $S_n$ denote the symmetric group of all permutations on
$\{1,\ldots,n\}$. For $\sigma\in S_n$ we let $\ve(\sigma)$
be $1$ or $-1$ according to whether $\sigma$ is even or odd.
For a matrix $A=(a_{ij})_{1\ls i,j\ls n}$ over a field the determinant
and the permanent of $A$ are defined by
$$\|A\|=\sum_{\sigma\in S_n}\ve(\sigma)\prod_{i=1}^na_{i,\sigma(i)}
\ \ \ \t{and}\ \ \ \per(A)=\sum_{\sigma\in
S_n}\prod_{i=1}^na_{i,\sigma(i)}$$
respectively.

Lemma 2.1 is very useful. For example, in view of Lemmas 1.1 and 2.1,
result (iii) follows from the following simple
observation:
$$\aligned\|x_j^{i-1}\|_{1\ls i,j\ls n}^*
=&\sum_{\sigma\in S_n}\ve(\sigma)\prod_{i=1}^n(x_{\sigma(i)})_{i-1}
\\=&\|(x_j)_{i-1}\|_{1\ls i,j\ls n}
=\|x_j^{i-1}\|_{1\ls i,j\ls n},
\endaligned\tag2.1$$ where in the last step we note that
$x^{r}\ (0\ls r<n)$ can be written as a linear combination of
$(x)_0,\ldots,(x)_r$.

Now we present our main technique
concerning the operator $P\mapsto P^*$.

\proclaim{Theorem 2.1}
Let $m_1,\ldots,m_n$ be nonnegative integers, and
let $A=(a_{ij})_{1\ls i,j\ls n}$ be a matrix over $\C$. Set
$$f(x_1,\ldots,x_n)=\|a_{ij}x_j^{m_i}\|_{1\ls i,j\ls n}P(x_1,\ldots,x_n)$$
where
$P(x_1,\ldots,x_n)\in\C[x_1,\ldots,x_n]$ is homogeneous and
$$P(x_1,\ldots,x_{i-1},x_j,x_{i+1},\ldots,x_{j-1},x_i,x_{j+1},\ldots,x_n)
=\nu P(x_1,\ldots,x_n)$$
for all $1\ls i<j\ls n$ with $\nu\in\{1,-1\}$.
Then
$$f^*(x,\ldots,x)=P^*(x-m_1,\ldots,x-m_n)\prod_{i=1}^n(x)_{m_i}
\times\cases\|A\|&\t{if}\ \nu=1,\\\per(A)&\t{if}\ \nu=-1.\endcases$$
\endproclaim
\Proof. Any $\sigma\in S_n$ can be written as a product of transpositions:
$$\sigma=(i_1j_1)\cdots(i_rj_r)\ \ \t{where}\ 1\ls i_s<j_s\ls n\
\t{for}\ s=1,\ldots,r.$$
Thus
$$P(x_{\sigma(1)},\ldots,x_{\sigma(n)})=\nu^rP(x_1,\ldots,x_n)
=\ve(\sigma)^{(1-\nu)/2}P(x_1,\ldots,x_n).$$

Write
$P(x_1,\ldots,x_n)=\sum_{j_1,\ldots,j_n}c_{j_1,\ldots,j_n}x_1^{j_1}\cdots x_n^{j_n}$.
Then
$$\align &f(x_1,\ldots,x_n)
\\=&\sum_{\sigma\in S_n}
\ve(\sigma)\prod_{i=1}^n\l(a_{i,\sigma(i)}x_{\sigma(i)}^{m_i}\r)\times P(x_1,\ldots,x_n)
\\=&\sum_{\sigma\in S_n}\(\ve(\sigma)\prod_{i=1}^n\l(a_{i,\sigma(i)}x_{\sigma(i)}^{m_i}\r)
\times\ve(\sigma)^{(\nu-1)/2}P(x_{\sigma(1)},\ldots,x_{\sigma(n)})\)
\\=&\sum_{\sigma\in S_n}\(\ve(\sigma)^{(\nu+1)/2}\prod_{i=1}^na_{i,\sigma(i)}
\times\sum_{j_1,\ldots,j_n}c_{j_1,\ldots,j_n}\prod_{i=1}^nx_{\sigma(i)}^{m_i+j_i}\).
\endalign$$
Therefore
$$\align f^*(x,\ldots,x)
=&\sum_{\sigma\in S_n}\ve(\sigma)^{(\nu+1)/2}\prod_{i=1}^na_{i,\sigma(i)}
\times\sum_{j_1,\ldots,j_n}c_{j_1,\ldots,j_n}\prod_{i=1}^n(x)_{m_i+j_i}
\\=&a\sum_{j_1,\ldots,j_n}c_{j_1,\ldots,j_n}\prod_{i=1}^n(x)_{m_i}
\times\prod_{i=1}^n(x-m_i)_{j_i}
\\=&a\prod_{i=1}^n(x)_{m_i}\times P^*(x-m_1,\ldots,x-m_n)
\endalign$$
where
$$a=\sum_{\sigma\in S_n}\ve(\sigma)^{(\nu+1)/2}\prod_{i=1}^na_{i,\sigma(i)}
=\cases\|A\|&\t{if}\ \nu=1,\\\per(A)&\t{if}\ \nu=-1.\endcases$$
This concludes the proof. \qed

\proclaim{Corollary 2.1} Let $m_1,\ldots,m_n$ be nonnegative integers.

{\rm (i) (Sun [Su])} If $A=(a_{ij})_{1\ls i,j\ls n}$ is a matrix with $a_{ij}\in\C$, then
$$\aligned&\(\|a_{ij}x_j^{m_i}\|_{1\ls i,j\ls n}
\prod_{1\ls i<j\ls n}(x_j-x_i)^{\da}\)^*(x,\ldots,x)
\\=&\prod_{1\ls i<j\ls n}(m_i-m_j)^{\da}\times
\prod_{i=1}^n(x)_{m_i}\times\cases\|A\|&\t{if}\ \da=0,\\\per(A)&\t{if}\ \da=1.\endcases
\endaligned\tag2.2$$

{\rm (ii) We have
$$\|x_j^{m_i}\|_{1\ls i,j\ls n}^*(x-n+1,\ldots,x)=\prod_{1\ls i<j\ls
n}(m_j-m_i)\times\prod_{i=1}^n\f{(x)_{m_i}}{(x)_{i-1}}.\tag2.3$$
\endproclaim
\Proof. As  $\|x_j^{i-1}\|_{1\ls i,j\ls n}=\prod_{1\ls i<j\ls n}(x_j-x_i)$ (Vandermonde),
by (2.1) we have
$$\align&\(\prod_{1\ls i<j\ls n}(x_j-x_i)^{\da}\)^*(x-m_1,\ldots,x-m_n)
\\=&\prod_{1\ls i<j\ls n}(x-m_j-(x-m_i))^{\da}=\prod_{1\ls i<j\ls n}(m_i-m_j)^{\da}.
\endalign$$
In view of this, Theorem 2.1 yields (2.2) immediately.

By Theorem 2.1,
$$\align&\l(\|x_j^{n-i}\|_{1\ls i,j\ls n}\times\|x_j^{m_i}\|_{1\ls i,j\ls n}\r)^*(x,\ldots,x)
\\=&n!\prod_{i=1}^n(x)_{n-i}\times\|x_j^{m_i}\|^*_{1\ls i,j\ls n}(x-n+1,\ldots,x).
\endalign$$
On the other hand, by part (i) we have
$$\align&\l(\|x_j^{n-i}\|_{1\ls i,j\ls n}\times\|x_j^{m_i}\|_{1\ls i,j\ls n}\r)^*(x,\ldots,x)
\\=&(-1)^{\bi n2}\(\|x_j^{m_i}\|_{1\ls i,j\ls n}
\times\prod_{1\ls i<j\ls n}(x_j-x_i)\)^*(x,\ldots,x)
\\=&n!\prod_{1\ls i<j\ls n}(m_j-m_i)\times\prod_{i=1}^n(x)_{m_i}.
\endalign$$
So (2.3) follows.

The proof of Corollary 2.1 is now complete. \qed

\Remark\ 2.1. When $m_i=(i-1)m$ for $i=1,\ldots,n$, Corollary 2.1(ii) yields
the following result related to [LS]:
$$\aligned&\(\prod_{1\ls i<j\ls n}(x_j^m-x_i^m)\)^*(x-n+1,\ldots,x)
\\=&1!2!\cdots (n-1)!m^{n(n-1)/2}\f{(x)_0(x)_m\cdots
(x)_{(n-1)m}}{(x)_0(x)_1\cdots(x)_{n-1}}.
\endaligned\tag2.4$$

\proclaim{Theorem 2.2} Let $m$ be any positive integer,
and let $a_1,\ldots,a_n$ be complex numbers.
Then
$$\aligned&\(\prod_{1\ls i<j\ls n}(x_j-x_i)^{2m-1}\)^*(x-n+1,\ldots,x)
\\=&(-1)^{(m-1)\bi
n2}\f{m!(2m)!\cdots(nm)!}{(m!)^nn!}
\times\f{(x)_0(x)_m\cdots(x)_{(n-1)m}}{(x)_0(x)_1\cdots(x)_{n-1}}
\endaligned\tag2.5$$
and
$$\aligned&\(\prod_{1\ls i<j\ls n}(a_jx_j-a_ix_i)(x_j-x_i)^{2m-1}\)^*(x,\ldots,x)
\\=&(-1)^{m\bi n2}\f{m!(2m)!\cdots(nm)!}{(m!)^nn!}
\per(a_j^{i-1})_{1\ls i,j\ls n}\prod_{r=0}^{n-1}(x)_{rm}.
\endaligned\tag2.6$$
\endproclaim
\Proof. Let $P_h(x_1,\ldots,x_n)=\prod_{1\ls i<j\ls n}(x_j-x_i)^h$
for $h=1,2,3,\ldots$. In light of Theorem 2.1,
$$\align P_{2m}^*(x,\ldots,x)=&\l((-1)^{\bi n2}\|x_j^{n-i}\|_{1\ls i,j\ls n}
P_{2m-1}(x_1,\ldots,x_n)\r)^*(x,\ldots,x)
\\=&(-1)^{\bi n2}n!\prod_{i=1}^n(x)_{n-i}\times P_{2m-1}^*(x-n+1,\ldots,x)
\endalign$$
and
$$\align&\l(\|a_j^{i-1}x_j^{i-1}\|_{1\ls i,j\ls n}P_{2m-1}(x_1,\ldots,x_n)\r)^*(x,\ldots,x)
\\=&\f{\per(a_j^{i-1})_{1\ls i,j\ls n}}{\per(a_j^0)_{1\ls i,j\ls n}}
\l(\|a_j^{0}x_j^{i-1}\|_{1\ls i,j\ls n}P_{2m-1}(x_1,\ldots,x_n)\r)^*(x,\ldots,x)
\\=&\f{\per(a_j^{i-1})_{1\ls i,j\ls n}}{n!}P_{2m}^*(x,\ldots,x).
\endalign$$
By Theorem 3.1 of Hou and Sun [HS],
$$P_{2m}^*(x,\ldots,x)=(-1)^{m\bi n2}\f{m!(2m)!\cdots(nm)!}{(m!)^n}(x)_0(x)_m\cdots(x)_{(n-1)m}.$$
So we have the desired (2.5) and (2.6). \qed

\proclaim{Corollary 2.2} Let $k,m,n$ be positive integers with
$k>m(n-1)$. Then
$$\aligned&[x_1^{k-n}\cdots x_n^{k-1}](x_1+\cdots+x_n)^{(k-1)n-mn(n-1)}
\prod_{1\ls i<j\ls n}(x_j-x_i)^{2m-1}
\\&=(-1)^{(m-1)\bi n2}\f{m!(2m)!\cdots(nm)!}{(m!)^nn!}
\cdot\f{((k-1-m(n-1))n)!}{\prod_{r=0}^{n-1}(k-1-rm)!}.
\endaligned\tag2.7$$
In particular,
$$\[\prod_{i=1}^nx_i^{(m-1)(n-1)+i-1}\]\prod_{1\ls i<j\ls n}(x_j-x_i)^{2m-1}
=(-1)^{(m-1)\bi n2}\f{(mn)!}{(m!)^nn!}.\tag2.8$$
\endproclaim

\Proof. Combining Lemma 2.1 with (2.5) we obtain (2.7).
(2.8) follows from (2.7) in the case $k=m(n-1)+1$. \qed

\Remark\ 2.3. Let $m_1,\ldots,m_n$ be nonnegative integers.
A confirmed conjecture of Dyson [D] can be stated as follows:
$$\align
&[x_1^{m_1(n-1)} \cdots x_n^{m_n(n-1)}]\prod_{1 \ls i < j \ls n} (x_i - x_j)^{m_i + m_j} \\
&\qquad\quad=(-1)^{\sum_{j=1}^n (j-1)m_j} \f{(m_1 + \cdots + m_n)!}{m_1!\cdots m_n!}.
\endalign$$
(See, e.g., Zeilberger [Z].)
Compared with this deep result, our (2.8) seems interesting too.

\heading{3. Proofs of Theorems 1.1--1.3}\endheading

\noindent{\it Proof of Theorem 1.1}. The case $n=1$ or
$k-1<m(n-1)$ is trivial. Below we assume $n\gs2$ and
$l=k-1-m(n-1)\gs0$.

As $|F|\gs p>mn\gs 2m$
we can extend each $S_{ij}\ (1\ls i<j\ls n)$ to a subset
$S_{ij}^*$ of $F$ with cardinality $2m-1$.
By Lemma 1.1 it suffices to show
that
$$[x_1^{k-n}\cdots x_n^{k-1}](x_1+\cdots+x_n)^{ln}
\prod_{1\ls i<j\ls n}\prod_{c\in S^*_{ij}}(x_j-x_i+c)$$
does not vanish. Let $e$ denote the multiplicative identity
of the field $F$. Then the above coefficient equals $he$ where
$$h=[x_1^{k-n}\cdots x_n^{k-1}](x_1+\cdots+x_n)^{ln}
\prod_{1\ls i<j\ls n}(x_j-x_i)^{2m-1}\in\Z.$$
By Corollary 2.2,
$$h=(-1)^{(m-1)\bi n2}\f{m!(2m)!\cdots(nm)!}{(m!)^nn!}
\cdot\f{(ln)!}{\prod_{r=0}^{n-1}(k-1-rm)!}.$$
As $p>mn$ and $p>ln$, $p$ does not divide $h$
and hence $he\not=0$. This concludes the proof. \qed

\medskip
\noindent{\it Proof of Theorem 1.2}. To avoid triviality,  we assume $n\gs 2$ and
$l=k-1-m(n-1)\gs0$.
As $q$ is odd,
the norms of those $1-\zeta_s/\zeta_t\ (1\ls s<t\ls n)$
(with respect to the field extension $\Q(e^{2\pi i/q})/\Q$)
are odd integers and hence $\|\zeta_t^{s-1}\|_{1\ls s,t\ls n}
=\prod_{1\ls s<t\ls n}(\zeta_t-\zeta_s)$
is not an algebraic integer times two.
Therefore $\per(\zeta_t^{s-1})_{1\ls s,t\ls n}\not=0$
as observed by Dasgupta et al. [DKSS]. By Lemma 2.1 and (2.6),
$$[x_1^{k-1}\cdots x_n^{k-1}](x_1+\cdots+x_n)^{ln}\prod_{1\ls s<t\ls n}
(\zeta_tx_t-\zeta_sx_s)(x_t-x_s)^{2m-1}\not=0.$$
Applying Lemma 1.1 we then obtain the desired
result. \qed

\medskip
\noindent{\it Proof of Theorem 1.3}.
 For $1\ls i<j\ls n$, let $r_{ij}$ denote
the unique integer in the interval $(-m_{ij}/2,m_{ij}/2]$
which is congruent to $b_i-b_j$ modulo $m_{ij}$.
For $x_i\in A_i$ and $x_j\in A_j$,
as $|x_i-x_j|<m_{ij}/2$ we have
 $$x_i+b_i\eq x_j+b_j\ (\mo\ m_{ij})
\iff x_j-x_i=r_{ij}.$$
Note also that
$$\per(\al_j^{i-1})_{1\ls i,j\ls n}
=\sum_{\sigma\in S_n}\prod_{i=1}^n\al_{\sigma(i)}^{i-1}>0.$$
Thus Theorem 1.3 follows from (2.6) and Lemmas 1.1 and 2.1. \qed
\medskip

\Ack. The main part of this work was done
during the first author's visit
to the second author's institute, Z. W. Sun
would like to thank the Institute of Mathematics, Academia Sinica
(Taiwan) for its support. The paper was revised during Sun's visit
to the University of California at Irvine, he is indebted to Prof.
Daqing Wan for the invitation.

\enddocument